\documentclass[12pt,leqno]{amsart}
\usepackage[dvips]{graphicx}
\usepackage{amsfonts}
\usepackage{amssymb}
\usepackage{amsmath}
\usepackage{amscd}
\usepackage{graphicx}

\def\DATE{\today}
\newtheorem{theorem}{Theorem}
\newtheorem{definition}[theorem]{Definition}

\newtheorem{lemma}[theorem]{Lemma}

\newtheorem{proposition}[theorem]{Proposition}

\newcommand\N{\mathbb{N}}
\newcommand\g{\mathfrak{g}}
\newcommand\h{\mathfrak{h}}

\newcommand\K{\mathbb{K}}

\newcommand\p{\mathcal{P}}

\newcommand\pf{\noindent{\it Proof. }}
\newcommand\lr{\left\{ \begin{array}{l}}
\email{michel.goze@uha.fr, elisabeth.remm@uha.fr}

\title{Lie algebras with associative structures. Applications to the study of $2$-step nilpotent Lie algebras}
\author{Michel Goze, Elisabeth Remm}
\date{}
\address{Universit\'{e} de Haute Alsace, LMIA, 4 rue des Fr\`{e}res Lumi\`{e}%
re, 68093 Mulhouse}
\email{Michel.goze@uha.fr , Elisabeth.Remm@uha.fr}
\begin{document}

\maketitle

\begin{abstract}
We investigate Lie algebras whose Lie bracket is also an associative or cubic associative multiplication to characterize the class of nilpotent Lie algebras with a  nilindex equal to $2$ or $3$.  In particular we study the class of $2$-step nilpotent Lie algebras, their deformations and we compute the cohomology which parametrize the deformations in this class.
\end{abstract}


\section{Introduction}
A finite dimensional Lie algebra $\g$ over a field of  characteristic zero  is $p$-nilpotent if for any $X \in \g$ the operator $ad X$ is $p$-nilpotent, that is, $(adX)^p=0. $ A very interesting class corresponds to $p=2,$ that is, the class of $2$-step nilpotent Lie algebras. This class is especially studied in the geometrical framework. In fact there are numerous studies on left invariant structures on a Lie group $G$ with an associated Lie algebra  $\g$ which is $2$-step nilpotent  (\cite{Eberlein}).  Here we are mostely interested by the algebraic study of the family  of $p$-step nilpotent Lie algebras for $p=2$ and $3.$ We show that such an algebra is defined by a Lie bracket which is also associative or $3$-associative. This leads to determine the properties of the corresponding operads. We show that the deformations of $2$-step nilpotent Lie algebras are governed by the operad cohomolgy and  we describe it. We compute this cohomology for two Lie algebras  $\mathfrak{k}_{2p+1}$ and $\mathfrak{k}_{2p}$ corresponding to the odd or even dimensional cases and  show that any  $2$-step nilpotent Lie algebra of dimension $2p+1$  (respectively $2p$) with maximal characteristic sequence is a linear deformation of $\mathfrak{k}_{2p+1}$ (respectively  $\mathfrak{k}_{2p}$). This permits a description of the class of $2$-step nilpotent Lie algebras with maximal characteristic sequence.

\section{Associative Lie multiplication}

Let $\g$ be a Lie algebra over a field of characteristic $0$. If we denote by $[X,Y]$ the Lie bracket of $\g$, it satisfies the following identities
$$
\left\{
\begin{array}{l}
\lbrack X,Y]=-[Y,X], \\
\lbrack \lbrack X,Y],Z]+\lbrack \lbrack Y,Z],X]+\lbrack \lbrack Z,X],Y]=0 \ \ {\text{ \rm (Jacobi Identity),}}
\end{array}
\right.
$$
for any $X,Y,Z \in \g$.
We assume moreover that the Lie bracket is also an associative product, that is, it satisfies
$$[[X,Y],Z]=[X,[Y,Z]],$$
for any $X,Y,Z \in \g$.
The Jacobi Identity therefore implies
$$[[Z,X ]],Y]= 0.$$
\begin{proposition}
The Lie bracket of the Lie algebra $\g$ is an associative product if and only if $\g$ is a two-step nilpotent Lie algebra.
\end{proposition}
In fact, the relation $[[Z,X ],Y]= 0$ means that $[[\g,\g],\g]=\mathcal{C}^2(\g)=0$ where $\mathcal{C}^i(\g)$ denotes the ideals of the descending central sequence of $\g$. The converse is obvious.

The classification of complex two-step nilpotent Lie algebras is known up to the dimension $7$. In the following our notation and terminology will be based on \cite{GK}

\smallskip

Any two-step nilpotent complex non abelian and indecomposable Lie algebra of dimension less than $7$ is isomorphic to one of the following algebras:

\begin{enumerate}

\item For the dimensions less than or equal to $3$:
\begin{itemize}
\item
$\frak{n}^1_3=\frak{h}_3$  \ : \ $\left[ X_1,X_2 \right]=X_3.$
\end{itemize}

\medskip

\item In dimension  $5$:
\begin{itemize}
\item $\frak{n}^5_5: \
\left[ X_1,X_2 \right]=X_3, \
\left[ X_1,X_4 \right]=X_5;$ \

\medskip

\item $\frak{n}^6_5= {\mbox {\rm the Heisenberg algebra $\frak{h}_2$}}: \
\left[ X_1,X_2 \right]=X_3, \
\left[ X_4,X_5 \right]=X_3.$ \

\end{itemize}

\medskip

\item In dimension $6$:
\begin{itemize}
\item
$\frak{n}^{19}_6: \
\left[ X_1,X_i \right]=X_{i+1}, \ i=2,4, \
\left[ X_2,X_6 \right]=X_5; \
$

\smallskip

\item $\frak{n}^{20}_6: \
\left[ X_1,X_i \right]=X_{i+1}, \ i=2,4, \
\left[ X_2,X_4 \right]=X_6. \
$
\end{itemize}

\medskip

\item In dimension $7$:

\begin{itemize}

\item $
\frak{n}_7^{120}:
\begin{array}{l}
\lbrack X_1, X_i]=X_{i+1}, \ i=2,4,6, \ \
\lbrack X_2 , X_4 ]= X_7;\
\end{array}
$

\smallskip

\item $
\frak{n}_7^{121}:
\begin{array}{l}
\lbrack X_1, X_i]=X_{i+1},  \ i=2,4,6; \
\end{array}
$

\smallskip

\item $
\frak{n}_7^{122}:
\begin{array}{l}
\lbrack X_1, X_i]=X_{i+1}, \ i=2,4,6, \ \
\lbrack X_4 , X_6 ]= X_7; \
\end{array}
$

\smallskip

\item $
\frak{n}_7^{123}:
\begin{array}{l}
\lbrack X_1, X_i]=X_{i+1}, \ i=2,4,6, \ \
\lbrack X_2, X_4 ]= X_5, \
\lbrack X_4 , X_6 ]= X_3;  \
\end{array}
$

\smallskip

\item $
\frak{n}_7^{124}:
\begin{array}{l}
\lbrack X_1, X_i]=X_{i+1}, \ i=2,4, \ \
\lbrack X_6, X_7 ]= X_5, \
\lbrack X_4 , X_7 ]= X_3; \
\end{array}
$

\smallskip

\item $
\frak{n}_7^{134} :
\begin{array}{l}
\lbrack X_1, X_i]=X_{i+1}, \ i=2,4, \ \
\lbrack X_6, X_7 ]= X_5; \
\end{array}
$

\medskip

\item$
\frak{n}_7^{126}:
\begin{array}{l}
\lbrack X_1, X_2]=X_{3},   \
\lbrack X_4 , X_5 ]= X_3, \
\lbrack X_6 , X_7 ]= X_3. \
\end{array}
$
\end{itemize}
\end{enumerate}
To develop  the operadic point of view, let us recall that an operad is a sequence $\p=\{\p(n), n \in \N^*\}$ of $\K[\Sigma _n]$-modules, where $\K[\Sigma _n]$ is the algebra group associated with the symmetric group $\Sigma_n$, with $comp_i$-operations (see  \cite{M.S.S}). The main example corresponds to the free operad $\Gamma(E)=\{\Gamma(E)(n)\}$ generated by a $\K[\Sigma_2]$-module. An operad $\p$ is called binary quadratic if there is a $\K[\Sigma_2]$-module $E$ and a $\K[\Sigma_3]$-submodule $R$ of $\Gamma(E)(3)$ such that $\p$ is isomorphic to $\Gamma(E)/\mathcal{R}$ where $\mathcal{R}$ is the operadic ideal generated by $\mathcal{R}(3)=R.$
\begin{proposition}There exists a binary quadratic operad, denoted by $2\mathcal{N}ilp$, with the property that any $2\mathcal{N}ilp$-algebra is a $2$-step nilpotent Lie algebra.
 \end{proposition}
 In fact, we consider $E=sgn_2,$ that is, the representation of $\Sigma_2$ by the signature, then $\Gamma(E)(3)=sgn_3 \oplus V_2$ where $V_2=\{(x,y,z) \in \K^3, x+y+z=0\}$. Let $R$ be the submodule of $\Gamma(E)(3)$ generated by
 $(x_i\cdot x_j)\cdot x_k$, $i,j,k$ all different. We deduce that $2\mathcal{N}ilp(2)$ is the $\K[\Sigma_2]$-module generated by $x_1\cdot x_2$ with the relation $x_2\cdot x_1=-x_1\cdot x_2$ and it is a $1$-dimensional vector space and $2\mathcal{N}ilp(3)=\{0\}$.
 \begin{proposition} \label{2Nilp}
The operad $2\mathcal{N}ilp$ is  Koszul.
\end{proposition}
Recall some general definitions and results on the duality of a binary quadratic operad (see \cite{M.S.S}).
The generating function of a binary quadratic operad $\p$
 is
$$g_{\p}(x)=\sum_{a \geq 1} \frac{1}{a!}\dim (\p(a))x^a.$$
Thus the generating function of the operad $2\mathcal{N}ilp$ is the polynomial
$$g_{2\mathcal{N}ilp}(x)=x+\frac{x^2}{2}.$$

The dual operad  $\p^!$ of the operad $\p$ is the quadratic operad $\p^! := \Gamma (E^\vee)/(R^\perp)$,
where $R^\perp \subset \Gamma (E^\vee)(3)$ is the annihilator of $R \subset \Gamma(E)(3)$ in the pairing
\begin{eqnarray}
\label{pairing}
\left\{
\begin{array}{l}
<(x_i \cdot x_j)\cdot x_k,(x_{i'} \cdot x_{j'})\cdot x_{k'}>=0, \ {\rm if} \ \{i,j,k\}\neq \{i',j',k'\}, \\
<(x_i \cdot x_j)\cdot x_k,(x_i \cdot x_j)\cdot x_k>=(-1)^{\varepsilon(\sigma)},  \\
\qquad \qquad \qquad   {\rm with} \ \sigma =
\left(
\begin{array}{lll}
i&j&k \\
i'&j'&k'
\end{array}
\right)
\ {\rm if} \ \{i,j,k\}= \{ i',j',k'\}, \\
<x_i \cdot (x_j\cdot x_k),x_{i'} \cdot (x_{j'}\cdot x_{k'})>=0, \ {\rm if} \ \{ i,j,k\} \neq \{ i',j',k'\}, \\
<x_i \cdot (x_j\cdot x_k),x_i \cdot (x_j\cdot x_k)>=-(-1)^{\varepsilon(\sigma)}, \\
\qquad \qquad \qquad  {\rm with} \ \sigma =
\left(
\begin{array}{lll}
i&j&k \\
i'&j'&k'
\end{array}
\right)
\ {\rm if} \ (i,j,k)=(i',j',k'), \\
<(x_i \cdot x_j)\cdot x_k,x_{i'} \cdot (x_{j'}\cdot x_{k'})>=0,
\end{array}
\right.
\end{eqnarray}
and $(R^\perp)$ is the
operadic ideal generated by $R^\perp.$ We deduce
$$\dim (2\mathcal{N}ilp)^!(1)=1, \ \dim (2\mathcal{N}ilp)^!(2)=1, \ \dim (2\mathcal{N}ilp)^!(3)=3, \ \dim (2\mathcal{N}ilp)^!(4)=15$$
and more generally, if we denote by $d_k$ the dimension of $(2\mathcal{N}ilp)^!(k)$, we have
$$
\left\{
\begin{array}{l}
\medskip
\displaystyle d_{2k+1}=\sum_{i=1}^k C_{2k+1}^i d_i d_{2k+1-i},\\
\displaystyle d_{2k}=\sum_{i=1}^{k-1} C_{2k}^i d_i d_{2k-i}+\frac{1}{2}C_{2k}^kd_k^2.
\end{array}
\right.
$$
In fact, the dual operad $(2\mathcal{N}ilp)^!$ is $\Gamma(1\! \! 1)$ the free operad generated by a commutative operation.
So the generating function of $2\mathcal{N}ilp^!$ is
$$\displaystyle \sum_{k \geq 1}\frac{d_k}{k!}x^k.$$
If an operad $\p$ is Koszul, then its dual $\p^!$ is also Koszul   and the generating functions are related by the functional equation
$$g_{\p}(-g_{\p ^!}(-x))=x.$$ It is known that $\Gamma(1\! \! 1)$ is Koszul, so also $2\mathcal{N}ilp $ and this implies the proposition.
We can verify that the generating function $g_{2\mathcal{N}ilp}$ of the operad $2\mathcal{N}ilp$ satisfies the  functional equation
$$g_{2\mathcal{N}ilp}(-g_{2\mathcal{N}ilp ^!}(-x))=x.$$

\medskip

\noindent{\bf Remarks.}
\begin{enumerate}
\item Recall that the operad satisfies the Koszul property if the corresponding free algebra is  Koszul, that is, its natural or operadic homology is trivial except in degree $0$. If $\mathcal{L}_r$ is the free Lie algebra of rank $r$
(i.e. on $k$ (free) generators), and if $\mathcal{C}^3(\mathcal{L}_r)$ is the third part of its descending central series, thus the free two-step nilpotent Lie algebra of rank $r$ is $\mathcal{N}(2,r)=\mathcal{L}_r/\mathcal{C}^3(\mathcal{L}_r)$. If $V_r$ is the $r$-dimensional vector space corresponding to the homogeneous component of degree $1$ of $\mathcal{L}_r$, thus $\mathcal{N}(2,r)=V_r \oplus  \bigwedge^2(V_r).$ For example, if $r=2$, thus $\mathcal{N}(2,2)$ is a $3$-dimensional Lie algebra with basis $\{e_1,e_2,e_1\wedge e_2\}$, where the last vector corresponds to $[e_1,e_2]$. If $r=3$, thus $\mathcal{N}(2,3)$ is the $6$-dimensional Lie algebra generated by $\{e_1,e_2,e_3,e_1\wedge e_2, e_1\wedge e_3,e_2\wedge e_3\}$. We can  define an homology of $\mathcal{N}(2,r)$ using the standard complex $(\bigwedge ^*(\mathcal{N}(2,r)), \partial_ *)$ where
$$\partial_p: \bigwedge ^p(\mathcal{N}(2,r))\rightarrow \bigwedge ^{p-1}(\mathcal{N}(2,r))$$
is defined by
$$\partial_p (x_1 \wedge x_2 \wedge \cdots \wedge x_p)= \displaystyle \sum_{i<j} (-1)^{i+j+1}[x_i,x_j]\wedge x_1 \wedge \cdots \widehat{x_i}\cdots \widehat{x_j}\cdots \wedge x_p$$
and $\partial_p =0$ is $p \leq 0.$ Let $m_p$ be the dimension of the $p$-th homology space 

\noindent $Ker \, \partial _p/ Im \, \partial_{p+1}$. For example if $r=2$, $\dim \mathcal{N}(2,2)=3$ and $\mathcal{N}(2,2)$ is generated by $e_1,e_2, [e_1,e_2]=e_3.$ It is isomorphic to the Heisenberg algebra. We have $m_0=1,m_1=m_2=2,m_3=1$. The general case was studied in \cite{Sigg}. The homology spaces are never trivial. For this complex, the free $2$-step nilpotent Lie algebra is not Koszul.
\item As a nilpotent Lie algebra is unimodular, we have the Poincar\'e duality. This implies that the second cohomology space of the free $2$-step nilpotent algebra of rank $r$
is trivial and this algebra is rigid in the variety of $2$-step nilpotent Lie algebra of dimension $r(r-1)/2$.  It defines an open orbit and an algebraic component in this variety.
\end{enumerate}

\medskip

\medskip

\noindent{\bf Remarks}
\begin{enumerate}
\item Let us consider an associative algebra $(A,\cdot)$ where $x \cdot y$ denotes the multiplication in $A$. Thus
$$[x,y]=x\cdot y - y\cdot x$$
is a Lie bracket. This Lie bracket is associative if and only if the multiplication of $A$ satisfies
$$(x \cdot y)\cdot z - (y \cdot x)\cdot z -(z \cdot x)\cdot y + (z \cdot y)\cdot x =0.$$
Let $v$ be the vector of $\K[\Sigma_3]$ $v=Id-\tau_{12}+\tau_{13}-c^2$ where $\tau_{ij}$ is the transposition  that exchanges the elements $i$ and $j$ and $c$ the $3$-cycle $(123)$. The orbit of $v$ with respect to the action of the group $\Sigma_3$ generates a $2$-dimensional vector space with basis  $\{v,\tau_{13}\cdot v\}$. We deduce that these algebras can be considered as $\p$-algebras where the quadratic operad $\p$ is defined by $\p(2)=sgn_2$ and $\p(3)=\Gamma(E)(3)/R$ with $R$ the submodule generated by $v((x_1x_2)x_3)$ and $\tau_{13}\cdot v((x_1x_2)x_3)$ with $\tau((x_1x_2)x_3)=((x_{\tau(1)}x_{\tau(2)})x_{\tau(3)})$ for any $\tau \in \Sigma_3$. In particular $\dim \p(3)=4.$
\item A Pre-Lie algebra is a non associative algebra defined by the identity
$$(xy)z-x(yz)=(xz)y-x(zy)$$
for all $x,y,z$. Assume that the Lie bracket of $\g$ satisfies also the Pre-Lie identity, that is,
$$[[x,y],z]-[x,[y,z]]=[[x,z],y]-[x,[z,y]].$$
Applying anticommutativity to this equation we obtain
$$[[x,y],z]+[[y,z],x]+[[z,x],y]+[[y,z],x]=0;$$
and finally the Jacobi identity gives
$$[[y,z],x]=0.$$
This shows that the Lie algebra is also $2$-step nilpotent and the Lie bracket is an associative product.
\item In \cite{Goze-Remm-Nonass}, we have defined classes of non associative algebras including in particular Pre-Lie algebras, Lie-admissible algebras and more generally  algebras with a non associative defining identity admitting a symmetry with respect to a subgroup of the symmetric group $\Sigma_3$. These algebras have been called $G_i$-associative algebras where $G_i$, $i=1, \cdots, 6$ are the subgroups of $\Sigma_3$. More precisely, a $G_1=\{Id\}$-associative algebra is an associative algebra, a $G_2=\{Id,\tau_{12}\}$-associative algebra is a Vinberg algebra that is, its multiplication  satisfies
    $$(xy)z-x(yz)=(yx)z-y(xz),$$
    a $G_3=\{Id,\tau_{23}\}$-associative algebra is a Pre-Lie algebra, a $G_4=\{Id,\tau_{13}\}$-associative algebra  satisfies
    $$(xy)z-x(yz)=(zy)x-z(yx),$$ a $G_5=\{Id,c,c^2\}$-associative algebra satisfies
    $$(xy)z-x(yz)+(yz)x-y(zx)+z(xy)-z(xy)=0,$$ and a $G_6=\Sigma_3$-associative algebra is a Lie-admissible algebra \cite{G.R1}. While writing this paper we discover that this notion already appear in \cite{Michor}.
    It is easy to see that if   the Lie bracket of $\g$ satisfies the
    $G_i$-associativity for $i=1,2,3,$ or $4$ then $\g$ is $2$-step nilpotent and the Lie bracket is an associative multiplication. The defining equations associated to the cases $i=5$ and $6$ are always satisfied because the $G_5$-conditions corresponds to the Jacobi identity and a Lie algebra is, in particular, a Lie-admissible algebra.
\end{enumerate}
\section{Cubic associative Lie multiplication}

Let $A$ be a $\K$ associative algebra with binary multiplication $xy$. The associativity which is the quadratic relation
$$(xy)z=x(yz)$$
implies six cubic relations
$$
\lr
((xy)z)t=(x(yz))t,\\
(x(yz))t=x((yz)t),\\
x((yz)t)=x(y(zt)), \ \ \ (*)\\
x(y(zt))=(xy)(zt),\\
(xy)(zt)=((xy)z)t.
\end{array}
\right.
$$
Recall that these relations correspond to the edges of the Stasheff pentagon.
\begin{definition}
\label{cubic associative}
A binary algebra, that is, an algebra whose multiplication is given by a bilinear map, is called cubic associative if the multiplication satisfies the cubic relations $(*)$.
\end{definition}
We call these relations cubic because if we denote by $\mu$ the multiplication, it occurs exactly three times in each  term of the relations. For example, the first relation writes as
$$\mu\circ(\mu \circ(\mu \otimes Id)\otimes Id)=\mu\circ(\mu \circ(Id \otimes \mu)\otimes Id)$$
which is cubic in $\mu$. It is the same thing for all other relations.

If $\mathcal{A}ss=\Gamma(E)/(R_{\mathcal{A}ss})$ is the operad for associative algebras, the relations $(*)$ are the generating relations of $(R_{\mathcal{A}ss})(4)$. But these relations are following from the relations defining $(R_{\mathcal{A}ss})(3)=R_{\mathcal{A}ss}$. In Definition \ref{cubic associative}, we do not assume that the algebra is associative. It is clear that $(*)$ do not implies associativity. From the relations $(*)$ we can define a binary cubic operad $\mathcal{A}ss\mathcal{C}ubic$. This operad will be studied in the next paragraph.

\begin{proposition}
Let $\g$ be a Lie algebra. The Lie bracket is cubic associative if and only if $\g$ is $3$-step nilpotent.
\end{proposition}
In fact, the first identity of $(*)$ becomes
$$
\begin{array}{lll}
[[[X_1,X_2],X_3],X_4] & = & [[X_1,[X_2,X_3]],X_4] \\
&=& -[[[X_2,X_3]],X_1],X_4]
\end{array}
$$
and finally
$$[[[X_1,X_2],X_3],X_4]+[[[X_2,X_3]],X_1],X_4]=[[[X_3,X_1]],X_2],X_4]=0,$$
which implies that $\g$ is $3$-nilpotent. Conversely, if $\g$ is $3$-nilpotent, all the relations of $(*)$ are satisfied.

The classification of the $3$-step nilpotent Lie algebra of dimension less than $7$ is the following
\medskip

\noindent{\bf Dimension $4$}

$\frak{n}^1_4$: $\left[ X_1,X_i \right]=X_{i+1}, \ i=2,3.$

\medskip

\noindent{\bf Dimension $5$}

\medskip

$\frak{n}^3_5$: $
\left[ X_1,X_i \right]=X_{i+1}, \ i=2,3, \
\left[ X_2,X_5 \right]=X_4;
$

\smallskip

$\frak{n}^4_5$: $
\left[ X_1,X_i \right]=X_{i+1}, \ i=2,3, \
\left[ X_2,X_3 \right]=X_5.
$

\medskip

\noindent{\bf Dimension $6$}

$\frak{n}^{11}_6$: $
\left[ X_1,X_i \right]=X_{i+1}, \ i=2,3,5, \
\left[ X_5,X_6 \right]=X_4; \
$

\smallskip

$\frak{n}^{12}_6$: $
\left[ X_1,X_i \right]=X_{i+1}, \ i=2,3,5, \
\left[ X_2,X_5 \right]=X_4; \
$

\smallskip

$\frak{n}^{13}_6$: $
\left[ X_1,X_i \right]=X_{i+1}, \ i=2,3,5, \
\left[ X_2,X_3 \right]=X_6, \
\left[ X_2,X_5 \right]=X_6; \
$

\smallskip

$\frak{n}^{14}_6$: $
\left[ X_1,X_i \right]=X_{i+1}, \ i=2,3,5, \
\left[ X_2,X_3 \right]=X_4-X_6, \
\left[ X_2,X_5 \right]=X_6; \
$

\smallskip

$\frak{n}^{15}_6$: $
\left[ X_1,X_i \right]=X_{i+1}, \ i=2,3,5, \
\left[ X_2,X_5 \right]=X_6, \
\left[ X_5,X_6 \right]=X_4; \
$

\smallskip

$\frak{n}^{16}_6$: $
\left[ X_1,X_i \right]=X_{i+1}, \ i=2,3,5, \
\left[ X_2,X_3 \right]=X_4; \
$

\smallskip

$\frak{n}^{17}_6$: $
\left[ X_1,X_i \right]=X_{i+1}, \ i=2,3,5; \
$
\medskip

$\frak{n}^{18}_6$: $
\left[ X_1,X_i \right]=X_{i+1}, \ i=2,3, \
\left[ X_5,X_6 \right]=X_4; \
$

\bigskip
\noindent{\bf Dimension $7$}

$
\frak{n}_7^{77} :
\begin{array}{l}
\lbrack X_1, X_i]=X_{i+1}, \ i=2,3,5,6, \ \
\lbrack X_2 , X_5 ]= X_7; \\
\end{array}
$

$
\frak{n}_7^{78} :
\begin{array}{l}
\lbrack X_1, X_i]=X_{i+1}, \ i=2,3,5,6, \ \
\lbrack X_2 , X_6 ]= X_4, \
\lbrack X_2 , X_5 ]= X_3; \\
\end{array}
$

\smallskip

$
\frak{n}_7^{79} :
\begin{array}{l}
\lbrack X_1, X_i]=X_{i+1}, \ i=2,3,5,6, \ \
\lbrack X_5 , X_6 ]= X_4, \
\lbrack X_2 , X_5 ]= X_7; \\
\end{array}
$

\smallskip

$
\frak{n}_7^{80} :
\begin{array}{l}
\lbrack X_1, X_i]=X_{i+1}, \ i=2,3,5,6; \ \\
\end{array}
$

\smallskip

$
\frak{n}_7^{81} :
\begin{array}{l}
\lbrack X_1, X_i]=X_{i+1}, \ i=2,3,5,6, \ \
\lbrack X_5 , X_6 ]= X_4; \\
\end{array}
$

\smallskip

$
\frak{n}_7^{82} :
\begin{array}{l}
\lbrack X_1, X_i]=X_{i+1}, \ i=2,3,5,6, \ \
\lbrack X_5 , X_6 ]= X_4, \
\lbrack X_2 , X_3 ]= X_7; \\
\end{array}
$

\smallskip

$
\frak{n}_7^{83} :
\begin{array}{l}
\lbrack X_1, X_i]=X_{i+1}, \ i=2,3,5,6, \ \
\lbrack X_5 , X_6 ]= X_7, \
\lbrack X_2 , X_3 ]= X_4+X_7; \\
\end{array}
$

\smallskip

$
\frak{n}_7^{84} :
\begin{array}{l}
\lbrack X_1, X_i]=X_{i+1}, \ i=2,3,5,6, \ \
\lbrack X_5 , X_6 ]= X_7, \
\lbrack X_2 , X_3 ]= X_4; \\
\end{array}
$

\smallskip

$
\frak{n}_7^{85} :
\left\{
\begin{array}{l}
\lbrack X_1, X_i]=X_{i+1}, \ i=2,3,5,6, \
\lbrack X_3 , X_5 ]= X_7, \
\lbrack X_2 , X_5 ]= X_4+X_6, \\
\lbrack X_2 , X_3 ]= X_4; \\
\end{array}
\right.
 $

\smallskip

$
\frak{n}_7^{86} :
\begin{array}{l}
\lbrack X_1, X_i]=X_{i+1}, \ i=2,3,5,6, \ \
\lbrack X_5 , X_6 ]= X_7, \
\lbrack X_2 , X_3 ]= X_7; \\
\end{array}
$

\smallskip

$
\frak{n}_7^{87} :
\left\{
\begin{array}{l}
\lbrack X_1, X_i]=X_{i+1}, \ i=2,3,5,6, \ \
\lbrack X_5 , X_6 ]= X_7+X_4, \
\lbrack X_2 , X_6 ]= X_4, \\
\lbrack X_2 , X_5 ]= X_3; \\
\end{array}
\right.
 $

\smallskip

$
\frak{n}_7^{88} :
\left\{
\begin{array}{l}
\lbrack X_1, X_i]=X_{i+1}, \ i=2,3,5,6, \ \
\lbrack X_5 , X_6 ]= X_4, \
\lbrack X_3 , X_5 ]= X_7, \\
\lbrack X_2 , X_3 ]= X_4,\
\lbrack X_2 , X_5 ]= X_6; \\
\end{array}
\right.
$

\smallskip

$
\frak{n}_7^{89} :
\left\{
\begin{array}{l}
\lbrack X_1, X_i]=X_{i+1}, \ i=2,3,5,6, \ \
\lbrack X_5 , X_6 ]= X_4, \
\lbrack X_2 , X_3 ]= X_4, \\
\lbrack X_2 , X_5 ]= X_7; \\
\end{array}
\right.
 $

\smallskip

$
\frak{n}_7^{90} :
\left\{
\begin{array}{l}
\lbrack X_1, X_i]=X_{i+1}, \ i=2,3,5,6, \ \
\lbrack X_2 , X_3 ]= X_4, \
\lbrack X_3 , X_5 ]= X_7, \\
\lbrack X_2 , X_5 ]= X_6; \\
\end{array}
\right.
 $

\smallskip

$
\frak{n}_7^{91} :
\begin{array}{l}
\lbrack X_1, X_i]=X_{i+1}, \ i=2,3,5,6, \ \
\lbrack X_5 , X_6 ]= X_7; \\
\end{array}
$

\smallskip

$
\frak{n}_7^{92} :
\begin{array}{l}
\lbrack X_1, X_i]=X_{i+1}, \ i=2,3,5,6, \ \
\lbrack X_2 , X_3 ]= X_4, \
\lbrack X_2 , X_5 ]= X_7; \\
\end{array}
$

\medskip

$
\frak{n}_7^{93} :
\begin{array}{l}
\lbrack X_1, X_i]=X_{i+1}, \ i=2,3,5, \ \
\lbrack X_2 , X_5 ]= X_7; \\
\end{array}
$

\smallskip

$
\frak{n}_7^{94} :
\begin{array}{l}
\lbrack X_1, X_i]=X_{i+1}, \ i=2,3,5, \ \
\lbrack X_2 , X_5 ]= X_4, \
\lbrack X_2 , X_3 ]= X_7; \\
\end{array}
$

\smallskip

$
\frak{n}_7^{95} :
\begin{array}{l}
\lbrack X_1, X_i]=X_{i+1}, \ i=2,3,5 \
\lbrack X_2 , X_3 ]= X_7; \\
\end{array}
$

\smallskip

$
\frak{n}_7^{96} :
\begin{array}{l}
\lbrack X_1, X_i]=X_{i+1}, \ i=2,3,5, \ \
\lbrack X_2 , X_5 ]= X_7, \
\lbrack X_2 , X_6 ]= X_4; \\
\end{array}
$

\smallskip

$
\frak{n}_7^{97} :
\begin{array}{l}
\lbrack X_1, X_i]=X_{i+1}, \ i=2,3,5, \ \
\lbrack X_2 , X_6 ]= X_4, \
\lbrack X_3 , X_5 ]= -X_4, \
\lbrack X_2 , X_5 ]= X_7; \\
\end{array}
$

\smallskip

$
\frak{n}_7^{98} :
\left\{
\begin{array}{l}
\lbrack X_1, X_i]=X_{i+1}, \ i=2,3,5, \ \
\lbrack X_2 , X_6 ]= X_4, \
\lbrack X_3 , X_5 ]= -X_4, \\
\lbrack X_2 , X_5 ]= X_7, \
\lbrack X_5 , X_6 ]= X_4; \\
\end{array}
\right.
$

\smallskip

$
\frak{n}_7^{99} :
\begin{array}{l}
\lbrack X_1, X_i]=X_{i+1}, \ i=2,3,5, \ \
\lbrack X_2 , X_7 ]= X_6, \
\lbrack X_2 , X_3 ]= X_4; \\
\end{array}
$

\smallskip

$
\frak{n}_7^{100} :
\begin{array}{l}
\lbrack X_1, X_i]=X_{i+1}, \ i=2,3,5, \ \
\lbrack X_2 , X_7 ]= X_4, \
\lbrack X_5 , X_7 ]= X_6; \\
\end{array}
$

\smallskip

$
\frak{n}_7^{101} :
\begin{array}{l}
\lbrack X_1, X_i]=X_{i+1}, \ i=2,3,5, \ \
\lbrack X_2 , X_7 ]= X_6, \
\lbrack X_5 , X_7 ]= X_4; \\
\end{array}
$

\smallskip

$
\frak{n}_7^{102} :
\begin{array}{l}
\lbrack X_1, X_i]=X_{i+1}, \ i=2,3,5, \ \
\lbrack X_5 , X_7 ]= X_4; \\
\end{array}
$

\smallskip

$
\frak{n}_7^{103} :
\begin{array}{l}
\lbrack X_1, X_i]=X_{i+1}, \ i=2,3,5, \ \
\lbrack X_2 , X_7 ]= X_4; \\
\end{array}
$

\smallskip

$
\frak{n}_7^{104} :
\begin{array}{l}
\lbrack X_1, X_i]=X_{i+1}, \ i=2,3,5, \ \
\lbrack X_5 , X_7 ]= X_4, \
\lbrack X_2 , X_3 ]= X_4; \\
\end{array}
$

\smallskip

$
\frak{n}_7^{105} :
\begin{array}{l}
\lbrack X_1, X_i]=X_{i+1}, \ i=2,3,5, \ \
\lbrack X_2 , X_7 ]= X_4, \
\lbrack X_2 , X_3 ]= X_4; \\
\end{array}
$

\smallskip

$
\frak{n}_7^{106} :
\begin{array}{l}
\lbrack X_1, X_i]=X_{i+1}, \ i=2,3,5, \ \
\lbrack X_2 , X_7 ]= X_4, \
\lbrack X_5 , X_6 ]= X_4; \\
\end{array}
$

\smallskip

$
\frak{n}_7^{107} :
\begin{array}{l}
\lbrack X_1, X_i]=X_{i+1}, \ i=2,3,5, \ \
\lbrack X_5 , X_7 ]= X_3, \
\lbrack X_6 , X_7 ]= X_4; \\
\end{array}
$

\smallskip

$
\frak{n}_7^{108} :
\begin{array}{l}
\lbrack X_1, X_i]=X_{i+1}, \ i=2,3,5, \ \
\lbrack X_2 , X_7 ]= X_6, \
\lbrack X_2 , X_3 ]= X_6; \\
\end{array}
$

\smallskip

$
\frak{n}_7^{109} :
\begin{array}{l}
\lbrack X_1, X_i]=X_{i+1}, \ i=2,3,5, \ \
\lbrack X_5 , X_7 ]= X_6, \
\lbrack X_2 , X_3 ]= X_6; \\
\end{array}
$

\smallskip

$
\frak{n}_7^{110} :
\begin{array}{l}
\lbrack X_1, X_i]=X_{i+1}, \ i=2,3,5, \ \
\lbrack X_2 , X_3 ]= X_6, \
\lbrack X_5 , X_7 ]= X_4;  \\
\end{array}
$

\smallskip

$
\frak{n}_7^{111} :
\begin{array}{l}
\lbrack X_1, X_i]=X_{i+1}, \ i=2,3,5, \ \
\lbrack X_2 , X_7 ]= X_6, \
\lbrack X_2 , X_5 ]= X_4; \\
\end{array}
$

\smallskip

$
\frak{n}_7^{112} :
\begin{array}{l}
\lbrack X_1, X_i]=X_{i+1}, \ i=2,3,5, \ \
\lbrack X_2 , X_3 ]= X_6, \
\lbrack X_2 , X_7 ]= X_4; \\
\end{array}
$

\smallskip

$
\frak{n}_7^{113} :
\begin{array}{l}
\lbrack X_1, X_i]=X_{i+1}, \ i=2,3,5, \ \
\lbrack X_5 , X_7 ]= X_6, \
\lbrack X_5 , X_6 ]= X_4; \\
\end{array}
$

\smallskip

$
\frak{n}_7^{114} :
\begin{array}{l}
\lbrack X_1, X_i]=X_{i+1}, \ i=2,3,5, \ \
\lbrack X_2 , X_7 ]= X_4, \
\lbrack X_5 , X_6 ]= X_4, \
\lbrack X_5 , X_7 ]= X_6; \\
\end{array}
$

\smallskip

$
\frak{n}_7^{115} :
\begin{array}{l}
\lbrack X_1, X_i]=X_{i+1}, \ i=2,3,5, \ \
\lbrack X_2 , X_5 ]= X_4, \
\lbrack X_5 , X_7 ]= X_3, \
\lbrack X_6 , X_7 ]= X_4; \\
\end{array}
 $

\smallskip

$
\frak{n}_7^{116} :
\left\{
\begin{array}{l}
\lbrack X_1, X_i]=X_{i+1}, \ i=2,3,5, \ \
\lbrack X_3 , X_5 ]= -X_4, \
\lbrack X_2 , X_6 ]= X_4,  \\
\lbrack X_5 , X_7 ]= -X_4; \\
\end{array}
\right.
$

\smallskip

$
\frak{n}_7^{117}(\alpha) :
\left\{
\begin{array}{l}
\lbrack X_1, X_i]=X_{i+1}, \ i=2,3,5, \ \
\lbrack X_2 , X_5 ]= X_7, \
\lbrack X_2 , X_7 ]= X_4, \\
\lbrack X_5 , X_6 ]= X_4, \
\lbrack X_5 , X_7 ]= \alpha X_4; \\
\end{array}
\right.
$

\smallskip

$
\frak{n}_7^{118} :
\left\{
\begin{array}{l}
\lbrack X_1, X_i]=X_{i+1}, \ i=2,3,5, \ \
\lbrack X_2 , X_5 ]= X_7, \
\lbrack X_2 , X_6 ]= X_4, \\
\lbrack X_3 , X_5 ]= -X_4, \
\lbrack X_5 , X_7 ]= -\frac{1}{4}X_4; \\
\end{array}
\right.
 $

\smallskip

$
\frak{n}_7^{119} :
\begin{array}{l}
\lbrack X_1, X_2]=X_3, \
\lbrack X_1 , X_3 ]= X_4, \
\lbrack X_2 , X_5 ]= X_4, \
\lbrack X_6 , X_7 ]= X_4. \\
\end{array}
$

\medskip

\noindent{\bf Remarks.}
\begin{enumerate}
\item We can generalize this process to define $(n-1)$-associative (binary) algebras: we consider the relations defining the $\Sigma_n$-module $\mathcal{A}ss(n)$ of the quadratic operad $\mathcal{A}ss$ and define, as above, an algebra with a multiplication which is a bilinear map $\mu$ (nonassociative), satisfying the previous relations where $\mu$ occurs $n-1$ times in each term of the relations. This algebra will be called $(n-1)$-associative (binary) algebra. If the Lie bracket of a algebra $\g$ is also $(n-1)$-associative, we prove a similar way than for the cubic associative case that $\g$ is a nilpotent Lie algebra of nilindex $n-1$.

\item There exits another notion of associativity for $n$-ary algebras (an $n$-ary algebra is a vector space with  a multiplication which is an $n$-linear map) , the total associativity. For example, a totally associative $3$-ary algebra has a ternary multiplication, denoted $xyz$, satisfying the relation:
    $$(xyz)tu=x(yzt)u=xy(ztu)$$
    for any $x,y,z,t,u$.
  The corresponding operad is studied in \cite{NG.R}, \cite{Nico-Elis-JOA}, \cite{Elis-tartu} and \cite{MR2}. Let $\g$ be a Lie algebra. We have the notion of Lie triple product given by $[[x,y],z].$ If we consider the vector space $\g$ provided with the $3$-ary product given by the Lie triple product, then $\g$ is a $3$-Lie algebra (\cite{Mic-Nic-Elis-Rabat} or \cite{Azcarraga}). Let us suppose now that the Lie triple bracket of $\g$ is a totally associative product. This implies
    $$
    [[[[X,Y],Z],T],U]=[[X,[[Y,Z],T]],U]=[[X,Y],[[Z,T],U]].$$
    But
    $$
    \begin{array}{lll}
    [[X,Y],[[Z,T],U]]&=&-[[[Z,T],U],[X,Y]]\\
    &=&[X,[Y,[[Z,T],U]]]+[Y,[[[Z,T],U],X]]\\
    &=& [[[[Z,T],U],Y],X]-[[[[Z,T],U],X],Y]\\
    &=&[[Z,[[T,U],Y]],X]-[[Z,[[T,U],X]],Y].
    \end{array}
    $$
    We deduce
    $$
    \begin{array}{lll}
    $$[[X,[[Y,Z],T]],U] &=&[[Z,[[T,U],Y]],X]-[[Z,[[T,U],X]],Y]\\
&=&2^5[[X,[[Y,Z],T]],U]-2^5[[X,[[Y,Z],U]],T]
    \end{array}
    $$
    Then
    $$(2^5-1)[[X,[[Y,Z],T]],U]=2^5[[X,[[Y,Z],U]],T].$$
    This implies
    $$[[X,[[Y,Z],T]],U]=0=[[[[X,Y],Z],T],U]=[[X,Y],[[Z,T],U]].$$
    The Lie algebra is $4$-step nilpotent.

\end{enumerate}

\section{Cubic operads}

Let $E$ be a $\K[\Sigma_2]$-module and $\Gamma(E)$ the free operad generated by $E$. Consider a $\K[\Sigma_4]$-submodule $R$ of $\Gamma(E)(4)$. Let $\mathcal{R}$ the ideal of $\Gamma(E)$ generated by $R$. We have
$$\mathcal{R}=\{\mathcal{R}(n), n \in \N^*\}$$
with $\mathcal{R}(1)=\{0\}, \ \mathcal{R}(2)=\{0\},  \ \mathcal{R}(3)=\{0\}, \ \mathcal{R}(4)=R.$
\begin{definition}
We call cubic operad generated by $E$ and defined by the relations $R\subset{\Gamma(E)(4)}$, the operad $\p(E,R)$ given by
$$\p(E,R)(n)=\displaystyle \frac{\Gamma(E)(n)}{\mathcal{R}(n)}.$$
\end{definition}
The operad $\mathcal{A}ss\mathcal{C}ubic$ is the cubic operad generated by $E=\K[\Sigma_2]$ and the $\K[\Sigma_4]$-submodule of relations $R$  generated by the vectors
$$
\left\{
\begin{array}{l}
((x_1x_2)x_3)x_4-(x_1(x_2x_3))x_4,
(x_1(x_2x_3))x_4-x_1((x_2x_3)x_4),
x_1((x_2x_3)x_4)-x_1(x_2(x_3x_4)),\\
x_1(x_2(x_3x_4))-(x_1x_2)(x_3x_4),
(x_1x_2)(x_3x_4)-((x_1x_2)x_3)x_4.
\end{array}
\right.
$$
Thus we have $\mathcal{A}ss\mathcal{C}ubic(2)=\K[\Sigma_2]$, $\mathcal{A}ss\mathcal{C}ubic(3)=\K[\Sigma_3]$, and
$$\mathcal{A}ss\mathcal{C}ubic(4)=\displaystyle \frac{\Gamma(E)(4)}{\mathcal{R}(4)}$$
is the $24$-dimensional $\K$-vector space generated by $\{((x_\sigma(1)x_\sigma(2))x_\sigma(3))x_\sigma(4), \ \sigma \in  \Sigma_4\}$.
\medskip

The operad $3\mathcal{N}ilp$ is the cubic operad defined by $3\mathcal{N}ilp(2)=sgn_2$, $3\mathcal{N}ilp(3)=sgn_3 \oplus V_2/sgn_3$ and $3\mathcal{N}ilp(4)=\{0\}$.
\begin{proposition}
The cubic operad $3\mathcal{N}ilp$ is not Koszul.
\end{proposition}
\pf In fact the Koszul operads has only quadratic relations.

\medskip

\noindent{\bf Remark: The Jordan operad.}  There is a cubic operad which is  really interesting, it is the operad $\mathcal{J}ord$ corresponding to Jordan algebras. Recall that a $\K$-Jordan algebra is a commutative algebra satisfying the following identity
$$x(yx^2)=(xy)x^2.$$
Since $\K$ is of zero characteristic, linearizing this identity, we obtain
$$((x_2x_3)x_4)x_1+((x_3x_1)x_4)x_2+((x_1x_2)x_4)x_3-(x_2x_3)(x_4x_1)-(x_3x_1)(x_4x_2)-(x_1x_2)(x_4x_3)=0.$$
This relation is cubic. It is invariant by the permutations $\tau_{12}, \tau_{13},\tau_{23}, c, c^2$ where $c$ is the cycle $(123).$ Thus the $\K[\Sigma_4]$-module $\mathcal{R}(4)$ generated by the vector
$$((x_2x_3)x_4)x_1+((x_3x_1)x_4)x_2+((x_1x_2)x_4)x_3-(x_2x_3)(x_4x_1)-(x_3x_1)(x_4x_2)-(x_1x_2)(x_4x_3)$$
is a vector space of dimension $4$. Let us consider the cubic operad $\mathcal{J}ord$ define by
$$\mathcal{J}ord(2)=1\! \! 1, \ \mathcal{J}ord(3)=1\! \! 1 \oplus V, \ \mathcal{J}ord(4)=\displaystyle \frac{\Gamma(1\! \! 1)(4)}{\mathcal{R}(4)}$$
where $1\! \! 1$ is the identity representation of $\Sigma_2$. Since $\dim \Gamma(1\! \! 1)(4)= 15$, thus
$\displaystyle \frac{\Gamma(1\! \! 1)(4)}{\mathcal{R}(4)}$ is of dimension $11$.

\medskip

\noindent If the Lie bracket of a Lie algebra $\g$ is also a Jordan product, then $\g$ is abelian. But maybe, it would be interesting to look Lie bracket satisfying the Jordan identity without the commutativity identity.
\section{The variety of $2$-step nilpotent Lie algebras}

\subsection{Chevalley-Hochschild cohomology of two-step nilpotent Lie algebras}

Let $\mu_0$ be a Lie bracket of a two-step nilpotent Lie algebra $\g_0$. If $\mu=\mu_0 + t \varphi$ is a linear deformation of $\mu_0$ then $\mu$ is a Lie bracket of a two-step nilpotent Lie algebra if and only if 
$\varphi$ satisfies the following conditions:

\begin{enumerate}
\item $\varphi \circ \varphi =0,$
\item $\delta_{H,\mu_0}(\varphi)=0,$ 
\item $\delta_{C,\mu_0}(\varphi)=0.$
\end{enumerate}
In the first condition, $\varphi$ is a skewsymmetric bilinear map and $\circ$ is the following trilinear map:
$$ \varphi \circ \varphi (X,Y,Z) =\varphi (\varphi (X,Y),Z).$$
In the second condition since the Lie algebra $\g_0$ is associative, $\delta_{H,\mu_0}$ is the coboundary
operator of the Hochschild cohomology (\cite{Schaps}), that is, 
$$\delta_{H,\mu_0}(\varphi)(X,Y,Z)=\mu_0(X, \varphi(Y,Z))-\varphi(\mu_0(X,Y),Z))+\varphi(X,\mu_0(Y,Z))-\mu_0(\varphi(X,Y),Z))$$
and in the third condition $\delta_{C,\mu_0}$ is the coboundary operator of the Chevalley operator (\cite{Ancochea}), that is, 
$$\begin{array}{rl}
\delta_{C,\mu_0}(\varphi)(X,Y,Z)= & \mu_0 (\varphi (X,Y),Z)+ \mu_0 (\varphi (Y,Z),X)+\mu_0 (\varphi (Z,X),Y) \\
& +\varphi (\mu_0 (X,Y),Z)+ \varphi (\mu_0 (Y,Z),X)+\varphi (\mu_0 (Z,X),Y).
\end{array}$$

\begin{lemma}
Let $\varphi$ be a skeysymmetric bilinear map on $\g_0$ then $\varphi$ satisfies the conditions (2) and (3) if and only if we have
$$\varphi(\mu_0(X,Y),Z)+\mu_0(\varphi(X,Y),Z)=0$$
for all $X,Y,Z$ in $\g_0.$
\end{lemma}

\pf In fact $\delta_{C,\mu_0}(\varphi)(X,Y,Z)+\delta_{H,\mu_0}(\varphi)(X,Y,Z)=\mu_0(\varphi(Z,X),Y)+\varphi(\mu_0(Z,X),Y)$ so $\delta_{C,\mu_0}(\varphi)(X,Y,Z)=\delta_{H,\mu_0}(\varphi)(X,Y,Z)=0$ implies $\varphi(\mu_0(X,Y),Z)+\mu_0(\varphi(X,Y),Z)=0.$ Conversely, considering
$$\delta_{H,C,\mu_0}(\varphi)(X,Y,Z)= \varphi(\mu_0(X,Y),Z)+\mu_0(\varphi(X,Y),Z)$$
then
$$\delta_{H,\mu_0}(\varphi)(X,Y,Z)=-\delta_{H,C,\mu_0}(\varphi)(X,Y,Z)-\delta_{H,C,\mu_0}(\varphi)(Y,Z,X)$$
and 
$$\delta_{C,\mu_0}(\varphi)(X,Y,Z)=\delta_{H,C,\mu_0}(\varphi)(X,Y,Z)+\delta_{H,C,\mu_0}(\varphi)(Y,Z,X)+\delta_{H,C,\mu_0}(\varphi)(Z,X,Y)$$
and $\delta_{H,C,\mu_0}(\varphi)=0$ implies that $\delta_{C,\mu_0}(\varphi)=\delta_{H,\mu_0}(\varphi)=0.$

In fact $\delta_{H,C,\mu_0}$ is a coboundary operator for the deformation cohomology of $2$-step nilpotent Lie algebras in the algebraic variety $2Nilp_n$ of $n$-dimensional $2$-step nilpotent Lie algebras.We assume now that $\mathbb{K}$ is algebraically closed. Recall that the variety $2Nilp_n$ is defined by the polynomial equations
$$\sum_{l=1}^n C_{ij}^lC_{lk}^s=0$$
for $1 \leq i,j,k,s\leq n,$ an element $(C_{ij}^k)$ with $1\leq i<j\leq n$ and $1 \leq k\leq n$ corresponds to a Lie algebra with constant structures $(C_{ij}^k)$ related to a fixed basis. Each element of this variety is an algebra on the operad $2\mathcal{N}ilp.$ From Proposition \ref{2Nilp} this operad is Koszul and the deformation cohomology is the classical operadic cohomology \cite{M.S.S}.

We deduce that this cohomology is associated with the complex $(\mathcal{C}^n(\mu_0,\mu_0), \delta^n_{H,C,\mu_0})$ where 
$\mathcal{C}^n(\mu_0,\mu_0)$ is the vector space of skew $n$-linear map on $\g_0$ with values in $\g_0$ and 
$\delta^n_{H,C,\mu_0}:\mathcal{C}^n(\mu_0,\mu_0) \rightarrow \mathcal{C}^{n+1}(\mu_0,\mu_0)$ is the linear operator defined by
$$\left\{
\begin{array}{ll}
\medskip
\delta^{2p}_{H,C,\mu_0}\psi(X_1, \cdots,X_{2p+1})= & \mu_0(X_1, \psi(X_2, \cdots,X_{2p+1}))\\
\medskip
& + \sum_{i=1}^p
\psi(X_1, \cdots,\mu_0(X_{2i},X_{2i+1}), \cdots ,X_{2p+1}),\\
\medskip
\delta^{2p-1}_{H,C,\mu_0}\psi(X_1, \cdots,X_{2p})= & \mu_0(X_1, \psi(X_2, \cdots,X_{2p}))\\
 & + \sum_{i=1}^{p-1}
\psi(X_1, \cdots,\mu_0(X_{2i+1},X_{2i+2}), \cdots ,X_{2p}).
\end{array}
\right.$$
In particular $\delta_{H,C,\mu_0}$ corresponds to $\delta^2_{H,C,\mu_0}.$ We will denote this cohomology by $H^*_{H,C}(\g_0,\g_0)$.

\subsection{Rigidity and deformations in $2Nilp_{2}$}

The algebraic linear group $Gl(n,\K)$ acts on the variety $2Nilp_n$, the orbit of an element $\mu_0$ associated with this action corresponds to the set of Lie algebras isomorphic to $\mu_0$ (we identify a Lie algebra with its Lie product). Thus, if we denote by $\mathcal{O}(\mu_0)$ this orbit, $\mu \in \mathcal{O}(\mu_0)$ if and only if there is $f \in Gl(n,\K)$ such that $\mu=f^{-1}\circ \mu_0 \circ (f \times f).$ 
\begin{definition}
A Lie algebra $\mu \in 2Nilp_n$ is called rigid in $2Nilp_n$ if its orbit $\mathcal{O}(\mu)$ is open in $2Nilp_n$ for the Zariski topology.
\end{definition}
\noindent{\bf Remarks}
\begin{itemize}
\item The rigidity in $2Nilp_n$ doesn't imply the rigidity in the algebraic variety $Nilp_n$ of all $n$-dimensional nilpotent $\K$-Lie algebra.
\item Today, we do not know nilpotent rigid Lie algebra in the full variety $Lie_n$ of $n$-dimensional Lie algebras. We do not know also  rigid Lie algebra in $Nilp_n$. It is a classical conjecture to affirm that no nilpotent Lie algebra is rigid in $Nilp_n$ or in $Lie_n$. 
\end{itemize}
The Nijenhuis-Richardson theorem writes in this case
\begin{theorem}
Let $\g_0$ be a $n$-dimensional $2$-step nilpotent Lie algebra on $\K$. If $H^2_{H,C}(\g_0,\g_0)=0$,  then $\g_0$ is rigid in $2Nilp_n$.
\end{theorem}
\begin{proposition}
The $(2p+1)$-dimensional Heisenberg algebra $\h_{2p+1}$ is rigid in $2Nilp_{2p+1}$.
\end{proposition}
\pf The dimension of the algebra of derivations of  $\h_{2p+1}$ is given, for example, in \cite{GP}. We have $\dim B^2_{H,C}(\h_{2p+1})=p(2p+1).$ Since the $2$-cochains are skew-symmetric, to compute 
$\delta_{H,C,\mu_0}(\varphi)=0$, it is sufficient to compute $\delta_{H,C,\mu_0}(\varphi)(X_i,X_j,X_k)=0$
for $i<j$, where $\{X_1,\cdots,X_{2p+1}\}$ is a basis of $\h_{2p+1}$ satisfying
$$[X_1,X_2]=\cdots=[X_{2i-1},X_{2i}]=\cdots=[X_{2p-1},X_{2p}]=X_{2p+1}.$$
If $\varphi(X_i,X_j)=\sum_{k=1}^{2p+1}a_{ij}^k X_k$,  then
$$\delta_{H,C,\mu_0}(\varphi)(X_i,X_j,X_k)=\varphi(\mu_0(X_i,X_j),X_k)+\mu_0(\varphi(X_i,X_j),X_k)=0$$is equivalent to
$$
\left\{
 \begin{array}{ll}
\medskip
\varphi(X_1,X_2)=\sum_{k=1}^{2p+1}a_{12}^k X_k;  \\
\medskip
\varphi(X_1,X_i)=a_{1i}^{2p+1} X_{2p+1}  , 3 \leq  i \leq 2p; \ \varphi(X_1,X_{2p+1})=-a_{12}^{1} X_{2p+1};\\
\medskip
\varphi(X_2,X_i)=a_{2i}^{2p+1} X_{2p+1}, \  3 \leq  i \leq 2p; \varphi(X_2,X_{2p+1})=a_{12}^{1} X_{2p+1}; \\
\cdots \\
\medskip
\varphi(X_{2i-1},X_{2i})=\sum_{k=1}^{2p}a_{12}^k X_k+a_{2i-1 \ 2i}^{2p+1}X_{2p+1}, \  2 \leq  i \leq p; \\
\medskip
 \varphi(X_l,X_{s})=a_{ls}^{2p+1} X_{2p+1}, (l,s) \neq (2i-1,2i);\\
\medskip
\varphi(X_{2l},X_{2p+1})=a_{12}^{2l-1} X_{2p+1}, \ l \leq p; \\
\medskip
\varphi(X_{2l-1},X_{2p+1})=-a_{12}^{2l} X_{2p+1} , \ l \leq p. 
\end{array}
\right.
$$
We deduce that $\dim Z^2_{H,C}(\h_{2p+1},\h_{2p+1})=p(2p+1).$
So $\dim H^2_{H,C}(\h_{2p+1},\h_{2p+1})=0$ and $\h_{2p+1}$ is rigid in $2Nilp_{2p+1}.$

\subsection{$2$-step nilpotent Lie algebras with characteristic sequence $(2,\cdots,2,1)$}

Recall that the characteristic sequence of a $n$-dimensional $\K$-Lie algebra $\g$ is  the invariant 
$$c(\g)=\displaystyle max_{X \in \g \backslash [\g,\g]}\{ (c_1(X),\cdots,c_p(X),1)\}$$
 where $(c_1(X),\cdots,c_p(X),1)$ is the decreasing sequence of the dimensions of the Jordan blocks of the nilpotent operator $ad(X)$, the maximum is computed with respect to the lexicographic order. This invariant as been introduced in \cite{Ancochea-Goze} to study classification of nilpotent Lie algebras (see also \cite{MGoran} and \cite{GK}). For $2$-step nilpotent Lie algebras this characteristic sequence is of type
$(2,\cdots,2,1,\cdots,1).$ For example the characteristic sequence of the Heisenberg algebra $\h_{2p+1}$ is $(2,1,\cdots,1).$ 
Let $\mathfrak{k}_{2p+1}$ be the $(2p+1)$-dimensional Lie algebras defined by the following brackets given in the basis $\{ X_1 ,\cdots , X_{2p+1} \}$ by
$$\left[X_1,X_{2i}\right]=X_{2i+1}, \ \ 1 \leq i \leq p,$$
the other non defined brackets are equal to zero. Its characteristic sequence is $(2,\cdots,2,1)$ where $2$ appears $p$ times.

\begin{lemma} \label{lemma1}
Any  $(2p+1)$-dimensional $2$-step nilpotent Lie algebra with characteristic sequence $(2,\cdots,2,1)$ is isomorphic to a linear deformation of $\mathfrak{k}_{2p+1}.$
\end{lemma}
 
\pf Let $\g$ be a $(2p+1)$-dimensional $2$-step nilpotent Lie algebra with characteristic sequence $(2,\cdots,2,1).$ There exists a basis $\{ X_1 ,\cdots ,X_{2p+1} \}$ such that the characteristic sequence is given by the operator $ad(X_1).$  If  $\{ X_1 ,\cdots ,X_{2p+1} \}$ is the Jordan basis of $ad(X_1)$ then the brackets of $\g$ write
$$\left\{
\begin{array}{l}
\medskip
[X_1,X_{2i}]=X_{2i+1}, \ \ 1 \leq i \leq p, \\
\left[X_{2i},X_{2j}\right]=\displaystyle\sum_{k=1}^pa_{2i,2j}^{2k+1}X_{2k+1}, \ \ 1\leq i<j \leq p.
\end{array}
\right.
$$
The change of basis $Y_1=X_1, \ Y_i=tX_i$ for $2\leq i \leq 2p+1$ shows that $\g$ is isomorphic to $\g_t$ whose brackets are 
$$\left\{
\begin{array}{l}
\medskip
[X_1,X_{2i}]=X_{2i+1}, \ \ 1 \leq i \leq p, \\
\left[X_{2i},X_{2j}\right]=\displaystyle t\sum_{k=1}^pa_{2i,2j}^{2k+1}X_{2k+1}, \ \ 1\leq i<j \leq p.
\end{array}
\right.
$$ 
If $\mu_t$ is the multiplication of $\g_t$ and $\mu_0$ the multiplication of $\mathfrak{k}_{2p+1}$ we have 
$\mu_t=\mu_0+t\varphi$ with $\varphi(X_{2i},X_{2j})=\sum_{k=1}^pa_{2i,2j}^{2k+1}X_{2k+1}, \ \ 1\leq i<j \leq p$ and $\varphi(X_l,X_s)=0$ for all the other cases with $l<s.$ So $\g$ is a linear deformation of $\mathfrak{k}_{2p+1}.$

Let us compute the second cohomological group $H_{H,C}^2(\mathfrak{k}_{2p+1},\mathfrak{k}_{2p+1}).$

Since for any $f $ in $End(\mathfrak{k}_{2p+1})$ we have 
$$\delta f(X_1,X_i)=-f(X_{2i+1})+a_{11}X_{2i+1}+
\displaystyle \sum_{k=1}^p a_{2k,2}X_{2k+1}$$ with $f(X_l)=\sum_{s=1}^{2p+1} a_{sl}X_s$, there exists in each class in $H_{H,C}^2(\mathfrak{k}_{2p+1},\mathfrak{k}_{2p+1})$ a representant $\varphi \in Z_{H,C}^2(\mathfrak{k}_{2p+1},\mathfrak{k}_{2p+1})$ such that $\varphi (X_1,  X_{2i})=0$ for $1\leq i
\leq p.$ This implies for such a cocycle that $0=\delta \varphi (X_1,X_{2i},X_1)=\varphi(X_{2i+1},X_1)$ and this cocycle satisfies $\varphi( X_1,Y)=0$ for any $Y$ in $\mathfrak{k}_{2p+1}.$ Moreover
$$0=\delta\varphi(X_1,X_{2i},X_l)=\varphi(X_{2i+1},X_l)$$ 
for any $1 \leq l \leq 2p+1$ and $ 1 \leq i \leq p.$ If we put $\varphi(X_{2i},X_{2j})=\sum_{k=1}^{2p+1} a_{2i,2j}^kX_k,$  the equations $\delta \varphi(X_{2i},X_{2j},X_l)=0$ for $l=1,2$ implies that $a_{2i,2j}^{2k}=a_{2i,2j}^{1}=0$ for $1\leq k\leq p.$ But $\delta f(X_{2i},X_{2j})=a_{1,2i}X_{2j+1}-a_{1,2j}X_{2i+1}.$ We can choose the cocycle $\varphi$ satisying 
$$\left\{
\begin{array}{l}
\medskip
\varphi (X_2,X_4)=\displaystyle \sum_{k=3}^{p}a_{24}^{2k+1}X_{2k+1} ,\\
\medskip
\varphi (X_2,X_{2i})=\displaystyle \sum_{k=2}^{p}a_{2,2i}^{2k+1}X_{2k+1}, \ 3\leq i\leq p, \\
\varphi (X_{2i},X_{2j})=\displaystyle \sum_{k=1}^{p}a_{2i,2j}^{2k+1}X_{2k+1} \ 2\leq i<j \leq p, 
\end{array}
\right.
$$
if $p \leq 3.$ If $p=2$ the chosen cocyle $\varphi$  is trivial and we obtain:

\begin{proposition}
Let $\mathfrak{k}_{2p+1}$ the $2$-step nilpotent Lie algebra defined by 
$[X_1,X_{2i}]=X_{2i+1}, \ \ 1 \leq i \leq p$ then 
\begin{itemize}
\item if $p=2,$ the algebra $\mathfrak{k}_{5}$ is rigid in $2Nilp_5,$
\item if $p>2,$ $\dim (H_{H,C}^2(\mathfrak{k}_{2p+1},\mathfrak{k}_{2p+1}))=\displaystyle \frac{p(p+1)(p-2)}{2}.$
\end{itemize}
\end{proposition}

We deduce that any $2$-step nilpotent $(2p+1)$-dimensional Lie algebra with characteristic sequence $(2,\cdots,2,1)$ is isomorphic to one of the following Lie algebras
$$\left\{
\begin{array}{l}
\medskip
\left[X_1,X_{2i}\right]=X_{2i+1}, \ 1 \leq i\leq p, \\
\medskip
\left[X_2,X_4\right]=\displaystyle \sum_{k=3}^{p}a_{24}^{2k+1}X_{2k+1}, \\
\medskip
\left[X_2,X_{2i}\right]=\displaystyle \sum_{k=2}^{p}a_{2,2i}^{2k+1}X_{2k+1}, \ 3\leq i\leq p ,\\
\left[X_{2i},X_{2j}\right]=\displaystyle \sum_{k=1}^{p}a_{2i,2j}^{2k+1}X_{2k+1} \ 2\leq i<j \leq p.
\end{array}
\right.
$$
Consider $\g$ a Lie algebra of this family $\mathcal{F}.$ Then the subspace $\frak{m}$  generated by $\left\{ X_2, \cdots, X_{2p+1}\right\}$ is a Lie subalgebra of $\g.$ So the classification up to an isomorphism  of the elements of the  family $\mathcal{F}$ corresponds to the classification of $2$-step nilpotent 
$(2p)$-dimensional  Lie algebras. Moreover we can assume that $X_2$ is a characteristic vector of $\frak{m}$ that is $c(\frak{m})$ is the characteristic sequence associated with $ad(X_2).$ But $[X_2,X_3]=0$ and $X_3 \notin Im  (ad\, X_2).$  We can assume that $\{ X_2, X_3,\cdots, X_{2p+1} \}$ is a Jordan basis of $ad\, X_2.$ For example, if $p=3$ we have $c(\frak{m})=(2,2,1,1)$ or $(2,1,1,1,1)$ or, in the abelian case, $(1,1,1,1,1,1,1).$  This corresponds to 
\begin{itemize}
\item $\left[X_2,X_{4}\right]=X_7, \left[X_2,X_{6}\right]=X_5$ if $c(\frak{m})=(2,2,1,1),$
\item $\left[X_2,X_{4}\right]=X_7$ if $c(\frak{m})=(2,1,1,1,1),$
\item  $\left[X_2,X_{i}\right]=0$ if $\frak{m}$ is abelian.
\end{itemize}

\subsection{$2$-step nilpotent Lie algebras with characteristic sequence $(2,\cdots,2,1,1)$}

Let us denote by $\frak{k}_{2p}$ the $(2p)$-dimensional Lie algebra given by the brackets
$$[X_1,X_{2i}]=X_{2i+1}, \ 1\leq i \leq p-1,$$
other non defined brackets are equal to zero.
\begin{lemma}
Any $2$-step nilpotent $(2p)$-dimensional Lie algebra is isomorphic to a linear deformation of $\frak{k}_{2p}.$
\end{lemma}

\pf It is similar to the proof of Lemma  \ref{lemma1}.

\medskip

If  $\mu_t=\mu_0+t\varphi$ is  a linear deformation of the bracket $\mu_0$ of  $\frak{k}_{2p}$ then 
$\varphi \in Z^2_{H,C}(\frak{k}_{2p},\frak{k}_{2p})$ and it is also a bracket of a $2$-step nilpotent Lie algebra. Let us determine these maps.

In $\frak{k}_{2p}$ we have 
$$\delta f(X_1,X_{2i})=-f(X_{2i+1})+a_{11}X_{2i+1}+\displaystyle \sum_{j=1}^{p-1}a_{2j,2i}X_{2j+1}$$
for $1 \leq i \leq p-1$ and 
$$\delta f(X_1,X_{2p})=\displaystyle \sum_{j=1}^{p-1}a_{2j,2p}X_{2j+1}.$$
Thus any $\tilde{\varphi} \in Z^2_{H,C}(\frak{k}_{2p},\frak{k}_{2p})$ is cohomologous to a cocycle $\varphi$ satisfying 
$$\left\{
\begin{array}{l}
\medskip
\varphi(X_1,X_{2i})=0, \ i=1,  \cdots, p-1, \\
\varphi(X_1,X_{2p})=\displaystyle \sum_{k=1}^{p}a_{1,2p}^{2k}X_{2k}.
\end{array}
\right.$$
For such a cocycle $\varphi$ we have 
$$0=\delta \varphi(X_1,X_{2p},X_1)=[\varphi(X_1,X_{2p}),X_1]$$
and $a^{2k}_{1, 2p}=0$ for $1 \leq k \leq p-1.$ This implies $\varphi(X_1,X_{2p})=a^{2p}_{1,2p}X_{2p}.$
Since $\mu_0+t\varphi $ is a multiplication of a $2$-step nilpotent Lie algebra, $a^{2p}_{1,2p}=0$ and  
$\varphi(X_1,X_{2i})=0$ for   $1 \leq i \leq p.$ This gives also 
$0=\delta \varphi(X_1,X_{2i},X_1)=\varphi(X_{2i+1},X_1)$ for $1 \leq i \leq p-1$ and $\varphi(X_1,Y)=0$ for any $Y \in \frak{k}_{2p}.$ This implies 
$0=\delta \varphi(X_1,X_{2i},X_j)=\varphi(X_{2i+1},X_{j})$ for $1 \leq i \leq p-1, \  1 \leq j \leq 2p.$
Thus this cocycle $\varphi$ satisfies
$$\varphi(X_{2i+1},Y)=0$$
for any $Y \in \frak{k}_{2p}$ and for any $1\leq i \leq p-1.$

We have also 
$0=\delta \varphi(X_{2i},X_{2j},X_1)=[\varphi(X_{2i},X_{2j}),X_1]$ and $\varphi(X_{2i},X_{2j})=\displaystyle \sum_{k=1}^{p-1} a_{2i,2j}^{2k+1}X_{2k+1} + a_{2i,2j}^{2p}X_{2p}.$ But
$$\left\{
\begin{array}{l}
\delta f(X_2,X_{2i})=a_{12}X_{2i+1}-a_{1, 2i} X_3, \ i=1,\cdots,p-1, \\
\delta f(X_2,X_{2p})=-a_{1,2p}X_3
\end{array}
\right.
$$
then we can assume that $\varphi$ satisfies
$$\left\{
\begin{array}{l}
\medskip
\varphi(X_1,X_{i})=0, \varphi(X_{2j+1},X_{k})=0,   i=1, \cdots, p, \ j=1,\cdots,p-1, \ k=1,\cdots ,2p, \\
\medskip
\varphi(X_2,X_{4})=\displaystyle \sum_{k=3}^{p-1}a_{2,4}^{2k+1}X_{2k+1}+a_{2,4}^{2p}X_{2p}, \\
\medskip
\varphi(X_2,X_{2i})=\displaystyle \sum_{k=2}^{p-1}a_{2,2i}^{2k+1}X_{2k+1}+a_{2,2i}^{2p}X_{2p}, \ i=3, \cdots, p \\
\varphi(X_{2i},X_{2j})=\displaystyle \sum_{k=1}^{p-1}a_{2i,2j}^{2k+1}X_{2k+1}+a_{2i,2j}^{2p}X_{2p}, \ 2\leq i<j \leq p.   
\end{array}
\right.
$$

By hypothesis the deformation $\mu_0+t\varphi$ is $2$-step nilpotent. This is equivalent to say that the cocycle $\varphi$ satisfies $\varphi (\varphi(X,Y),Z)=0$ for any $X,Y,Z.$ This gives
$$a_{2i,2j}^{2p}\varphi(X_{2k},X_{2p})=0$$ for any $1 \leq i,j,k \leq p.$
In particular $a_{2i,2p}^{2p}\varphi(X_{2i},X_{2p})=0$  which implies that $a_{2i,2p}^{2p}=0$ for any $1 \leq i \leq p.$ If there exits $(i,j)$ with $a_{2i,2j}^{2p}\neq 0$ then $\varphi (X_{2k},X_{2p})=0$ for any $1 \leq k \leq p.$ In this case the Lie algebra $\g$ is defined by the cocycle 
$$\left\{
\begin{array}{l}
\medskip
\varphi(X_1,X_{i})=0, \varphi(X_{2j+1},X_{k})=0,   i=1, \cdots, p, \ j=1,\cdots,p-1, \ k=1,\cdots ,2p, \\
\medskip
\varphi(X_2,X_{4})=\displaystyle \sum_{k=3}^{p-1}a_{2,4}^{2k+1}X_{2k+1}+a_{2,4}^{2p}X_{2p}, \\
\medskip
\varphi(X_2,X_{2i})=\displaystyle \sum_{k=2}^{p-1}a_{2,2i}^{2k+1}X_{2k+1}+a_{2,2i}^{2p}X_{2p}, \ i=3, \cdots, p-1 \\
\varphi(X_{2i},X_{2j})=\displaystyle \sum_{k=1}^{p-1}a_{2i,2j}^{2k+1}X_{2k+1}+a_{2i,2j}^{2p}X_{2p}, \ 2\leq i<j \leq p-1.   
\end{array}
\right.
$$
If $a_{2i,2j}^{2p}=0$ for any $1 \leq i,j \leq p$ then $\g$ is defined by the cocycle 
$$\left\{
\begin{array}{l}
\medskip
\varphi(X_1,X_{i})=0, \varphi(X_{2j+1},X_{k})=0,   i=1, \cdots, p, \ j=1,\cdots,p-1, \ k=1,\cdots ,2p, \\
\medskip
\varphi(X_2,X_{4})=\displaystyle \sum_{k=3}^{p-1}a_{2,4}^{2k+1}X_{2k+1}, \\
\medskip
\varphi(X_2,X_{2i})=\displaystyle \sum_{k=2}^{p-1}a_{2,2i}^{2k+1}X_{2k+1}, \ i=3, \cdots, p \\
\varphi(X_{2i},X_{2j})=\displaystyle \sum_{k=1}^{p-1}a_{2i,2j}^{2k+1}X_{2k+1}, \ 2\leq i<j \leq p.   
\end{array}
\right.
$$
In particular thiese Lie algebras have been classified in dimension $8$ in   \cite{Zaili}.

\end{document}